\documentclass[12pt,reqno]{amsart}
\usepackage{amsfonts,amsmath,amsthm,amssymb}
\usepackage[usenames,dvipsnames]{color}
\usepackage{enumerate,nicefrac}
\usepackage{hyperref,multicol}
\usepackage{rotating}
\usepackage[labelsep=space]{caption}
\usepackage[all]{xy}
\usepackage{tikz}
\usetikzlibrary{decorations.markings}
\usepackage[labelformat=empty]{caption}
\usepackage{mathtools}

\newtheorem{theorem}{Theorem}[section]

\newtheorem{proposition}[theorem]{Proposition}

\theoremstyle{definition}

\usepackage[ansinew]{inputenc}
\usepackage{ifthen}
\usepackage{animate}

\newcommand{\be}{\begin{equation}}
\newcommand{\ee}{\end{equation}}
\newcommand{\NZ}{\reflectbox{Z}}

\title[Representations of imaginary quadratic integers]{On the representation of an imaginary quadratic integer in two different bases}

\author{Divyum Sharma}
\address{Department of Mathematics\\
Birla~Inst{i}tute~of~Technology~and~Science, Pilani 333\,031 \textsc{India}}
\email{divyum.sharma\symbol{64}pilani.bits-pilani.ac.in}

\begin{document}

\date{}

\begin{abstract}
Let $(\alpha,\mathcal{N}_{\alpha})$ and $(\beta,\mathcal{N}_{\beta})$ be two canonical number systems for an imaginary quadratic number field $K$ such that $\alpha$ and $\beta$ are multiplicatively independent. We provide an effective lower bound for the sum of the number of non-zero digits in the $\alpha$-adic and $\beta$-adic expansions of an algebraic integer $\gamma\in\mathcal{O}_K$ which is an increasing function of $|\gamma|$. This is an analogue of an earlier result due to Stewart on integer representations.
\end{abstract}

\maketitle
\section{Introduction}
\noindent In 1971, Senge and Straus \cite{SeSt1} (also in \cite{SeSt2}) proved that the number of integers with bounded sum of digits with respect to two different bases $q_1$ and $q_2$ is finite if and only if $\log q_1/\log q_2$ is irrational. Stewart \cite{St_80} obtained an effective version of this result  using Baker's inequality for linear forms in logarithms of algebraic numbers, and a result of Loxton and van der Poorten on multiplicative relations in number fields. As remarked by Bugeaud, Cipu and Mignotte \cite{BCM_13}, this result is an illustration of the principle that \textit{in two unrelated number systems, two miracles cannot happen simultaneously for large integers}. We refer to \cite{BHLS,EGST,Lu_00,Mi88,Sc90,Sp23} and the references therein for further results on this theme. In this paper, we explore the same problem for number systems in the ring of algebraic integers of an imaginary quadratic field.
\par Let $K$ be an algebraic number field of degree $n$ and $\mathcal{O}_K$ be its ring of integers.  Further, let $\alpha\in\mathcal{O}_K \setminus \{ 0 \}$ and $\mathcal{N}_{\alpha} := \{\,0, 1, \ldots, |N_{K/\mathbb{Q}}(\alpha)|-1\,\}$. The pair $(\alpha,\mathcal{N}_{\alpha})$ is said to be a \textit{canonical number system} (CNS) for $\mathcal{O}_K$ if each $\gamma\in\mathcal{O}_K$ admits a unique representation of the form
      \[
        \gamma=\sum\limits_{0\leq i\leq L} \epsilon_i \alpha^i,
     \]
 where $L\geq 0$, $\epsilon_{L} \neq 0$ if $L \neq 0$ and $\epsilon_i \in \mathcal{N}_{\alpha}$ for $0 \leq i \leq L$.  
 
 When $K=\mathbb{Q}(i)$, K\'{a}tai and Szab\'{o} \cite{KS_75}  proved that $(\alpha,\mathcal{N}_{\alpha})$ is a CNS if and only if Re($\alpha$) is a negative integer and Im$(\alpha)=\pm 1$.
More generally, if $[K:\mathbb{Q}]=2$ and the minimal polynomial of $\alpha\in\mathcal{O}_K$ is $x^2+Ex+F$, Gilbert \cite{Gi_81} and independently K\'{a}tai and  Kov\'{a}cs \cite{KK_80,KK_81} proved that $(\alpha,\mathcal{N}_{\alpha})$ is a CNS if and only if
    \[
    F\geq 2\ \textrm{ and } -1\leq E\leq F.
    \]
 We refer to \cite{Mo_10,Mo_12,Th_98,Th_00} and the references therein for results related to the distribution of the sum-of-digits functions of CNS in number fields.
 For such an $\alpha$ and for $\gamma\in \mathcal{O}_K$, let $S(\gamma;\alpha)$ and $\NZ(\gamma;\alpha)$ denote  the sum of digits and the number of non-zero digits in the $\alpha$-adic representation of $\gamma$, respectively.
 Note that $\NZ(\gamma;\alpha)\leq S(\gamma;\alpha)$ for all $\gamma\in\mathcal{O}_K$.
 In this note, we prove the following result.
\begin{theorem}\label{thm1}
    Let $K$ be an imaginary quadratic number field. Further, let $(\alpha,\mathcal{N}_{\alpha})$ and $(\beta,\mathcal{N}_{\beta})$ be two canonical number systems in $K$. If $\alpha$ and $\beta$ are multiplicatively independent, then for all $\gamma$ in $\mathcal{O}_K$, 
    \begin{enumerate}[(i)]
     \item $\NZ(\gamma;\alpha)+\NZ(\gamma;\beta)\geq \frac{\log\log |\gamma|}{\log\log\log |\gamma|+C}$, and hence
        \item $S(\gamma;\alpha)+S(\gamma;\beta)\geq \frac{\log\log |\gamma|}{\log\log\log |\gamma|+C}$,
    \end{enumerate}
    where $C$ is an effectively computable positive number depending only on $\alpha$ and $\beta$.
\end{theorem}
\noindent Recall that the complex numbers $z_1,z_2$ are said to be \textit{multiplicatively dependent} if there exist integers $u_1,u_2$, not both zero, such that $$z_1^{u_1}z_2^{u_2}=1.$$ Otherwise, we say that $z_1,z_2$ are \textit{multiplicatively independent}.

\textbf{Remark. }
\begin{enumerate}
    \item If $\alpha$ and $\beta$ are multiplicatively dependent with $\alpha^u=\beta^v$ (say) where $u,v\in\mathbb{N}$, then the above result does not hold as for every positive integer $m$, we have
    \[
   \NZ(\alpha^{um};\alpha)=1=\NZ(\alpha^{um};\beta),
    \]
   where the last equality holds because $\alpha^{um}=\beta^{vm}$. 
    \item As mentioned earlier, $(\alpha,\mathcal{N}_{\alpha})$ is a CNS in $\mathbb{Q}(i)$ if and only if $\alpha=-a\pm i$ for some positive integer $a$. If $\alpha$ and $\beta$ are multiplicatively dependent, with
    \[
    \alpha=-a\pm i,\beta=-b\pm i, \textrm{ for some } a,b\in\mathbb{N},\  a\neq b,
    \]
    then the integers $a^2+1,b^2+1$ are multiplicatively dependent. Hence $a^2+1=x^v$ and $b^2+1=x^w$ for some rational integers $x,v,w$ with $x>0$ and $\max(v,w)\geq 3$. Using Lebesgue's theorem \cite{Le1850}, it follows that this is not possible. Hence the bases $\alpha=-a\pm i$ and $\beta=-b\pm i$ are   multiplicatively independent if $a\neq b$. Indeed, it can be seen that the only multiplicatively dependent bases for  CNS in $\mathbb{Q}(i)$  are $-1+i$ and $-1-i$.
\end{enumerate}

In Section \ref{prelim}, we record certain preliminary results required towards our proof of Theorem \ref{thm1}. In Section \ref{proof}, we complete the proof following Stewart's strategy. The idea is to use a Baker-type result on linear forms in logarithms of algebraic numbers along with a result of Loxton and van der Poorten on multiplicative relations in number fields.  On certain occasions, we make use of a result due to Kov\'{a}cs and Peth\H{o} \cite{KP_92} (Proposition \ref{highestexp} below)  about the length of
 $\alpha$-adic expansions for algebraic integers. This is to obtain certain inequalities whose analogues were readily available in the case of rational integers due to the ordering relation therein. 
\section{Preliminaries}\label{prelim}

We first record a result due to Kov\'{a}cs and Peth\H{o} \cite{KP_92} about the length of
the $\alpha$-adic expansion of an element $\gamma \in\mathcal{O}_K$.
\begin{proposition}\label{highestexp}
    Let $(\alpha,\mathcal{N}_{\alpha})$ be a canonical number system in a number field $K$ of degree $D$ with ring of integers $\mathcal{O}_K$. Then there exist effectively computable constants $e_1(\alpha),e_2(\alpha)$, depending only on $\alpha$ such that if 
    $\gamma=\sum_{i=0}^L\epsilon_i\alpha^i$ is the representation of a non-zero algebraic integer $\gamma$ in the number system $(\alpha,\mathcal{N}_{\alpha})$, then
    \[
    \max_{1\leq i\leq D}\frac{\log|\gamma^{(i)}|}{\log|\alpha^{(i)}|}+e_1(\alpha)\leq L\leq \max_{1\leq i\leq D}\frac{\log|\gamma^{(i)}|}{\log|\alpha^{(i)}|}+e_2(\alpha),
    \]
    where $\alpha^{(1)},\ldots,\alpha^{(D)}$ and $\gamma^{(1)},\ldots,\gamma^{(D)}$ denote the conjugates of $\alpha$ and $\gamma$, respectively.
\end{proposition}
If $K$ is an imaginary quadratic field, the above inequality may be written as
    \begin{equation}\label{length}
    \begin{aligned}
          \frac{\log|\gamma|}{\log|\alpha|}+e_1(\alpha)&\leq L\leq \frac{\log|\gamma|}{\log|\alpha|}+e_2(\alpha),\\
          \textrm{or, }\  |\alpha|^{L-e_2(\alpha)}&\leq |\gamma|\leq |\alpha|^{L-e_1(\alpha)}.
    \end{aligned}
         \end{equation}
(Note that it follows from the result of K\'{a}tai and  Kov\'{a}cs \cite{KK_80,KK_81} and Gilbert \cite{Gi_81} mentioned earlier that if $(\alpha,\mathcal{N}_{\alpha})$ is a CNS in an imaginary quadratic field, then $|\alpha|>1$.)

Next, we recall the definition and some basic properties of the absolute logarithmic height function $h(\cdot)$. If $\delta$ is an algebraic  number with minimal polynomial $f(X)=d_0(X-\delta^{(1)})\cdots(X-\delta^{(d)})\in\mathbb{Z}[X]$, 
    \[
    h(\delta)=\frac{1}{d}\left(\log d_0+\sum_{i=1}^d\max(0,\log|\delta^{(i)}|)\right)
    \]
denotes its absolute logarithmic height.
If $\delta_1,\delta_2$ are algebraic numbers, then
    \begin{align}\label{height_of_product}
        h(\delta_1\delta_2^{\pm 1})&\leq h(\delta_1)+h(\delta_2)
    \end{align}
(see \cite[Property 3.3]{Wa_00}). The next result gives an upper bound for the absolute logarithmic height of an algebraic number which is given as the value of a polynomial in an algebraic number.
\begin{proposition}\cite[Special case of Lemma 3.7]{Wa_00}\label{height_of_poly}
    Let $f(x)\in\mathbb{Z}[x]$ be a non-zero polynomial and let $\delta$ be an algebraic number. Then
    \[
    h(f(\delta))\leq (\deg f)h(\delta)+\log L(f),
    \]
    where $L(f)$ denotes the sum of the absolute values of the coefficients of $f$.
\end{proposition}
To prove our results, we will use a Baker-type lower bound for a non-zero linear form in logarithms of algebraic numbers. Following is a general estimate due to Matveev \cite{Mat1, Mat2}. 
\begin{proposition}\label{lfl}
    Let $T$ be a positive integer. Let  $\delta_1,\ldots,\delta_T$ be non-zero algebraic numbers and $\log\delta_1,\ldots,\log\delta_T$ be some determinations of their complex logarithms. Let $D$ be the degree of the number field generated by $\delta_1,\ldots,\delta_T$ over $\mathbb{Q}$.  For $j=1,\ldots,T$, let $A_j'$ be a real number satisfying
    \[
    \log A_j'\geq\max\left(h(\delta_j),\frac{|\log\delta_j|}{D},\frac{0.16}{D}\right).
    \]
    Let $k_1,\ldots,k_T$ be rational integers. Set
    \[
    B=\max(|k_1|,\ldots,|k_T|) \textrm{ and }\ \Lambda=\delta_1^{k_1}\cdots\delta_T^{k_T}-1.
    \]
    Then
    \[
    \log|\Lambda|>-3\times 30^{T+4}(T+1)^{5.5}D^{T+2}\log(eD)\log A_1'\cdots\log A_T'\log(eTB).
    \]
\end{proposition}
To deal with the degenerate case when $\Lambda=0$, we will use the following result about \textit{`small'} multiplicative relations between multiplicatively dependent algebraic numbers due to Loxton and van der Poorten \cite{LV_83}.
\begin{proposition}\label{multdep}
  Let $T\geq 2$ be an integer and $\delta_1,\ldots,\delta_T$ be algebraic numbers in an algebraic number field $K$ of degree $D$. Let $\omega(K)$ denote the number of roots of unity in $K$.
Suppose there are rational integers $u_1,\ldots,u_T$, not all zero, such that 
\[
\delta_1^{u_1}\cdots\delta_T^{u_T}=1.
\] 
Then there exist rational integers $n_1,\ldots,n_T$, not all zero, such that 
    \[
    \delta_1^{n_1}\cdots\delta_T^{n_T}=1
    \]
and
     \[
    |n_j|\leq (T-1)!\omega(K)\prod_{i\neq j}(Dh(\delta_i)/\lambda(D)),\ j=1,\ldots,T,
    \]
    where $\lambda(D)$ is a positive number depending only on $D$.
\end{proposition}

\section{Proof of Theorem \ref{thm1}}\label{proof}
Define the complex logarithm by $$\log z=\log |z|+i\,\arg z,\ -\pi<\arg z\leq \pi.$$
Without loss of generality, assume that $|\alpha|\geq |\beta|$. 
Let $\gamma\in\mathcal{O}_K$. Write 
    \[    \gamma=a_{m_1}\alpha^{m_1}+a_{m_2}\alpha^{m_2}+\cdots+a_{m_r}\alpha^{m_r},
    \]
    where $m_1>m_2>\cdots>m_r\geq 0,\ a_i\in\mathcal{N}_{\alpha}\setminus\{0\}\ \textrm{ for }\ i=1,\ldots, r$.
    In particular, $\NZ(\gamma;\alpha)=r.$
Similarly,
    \[    \gamma=b_{l_1}\beta^{l_1}+b_{l_2}\beta^{l_2}+\cdots+b_{l_t}\beta^{l_t},
    \]
    where $l_1>l_2>\cdots>l_t\geq 0,\ b_i\in\mathcal{N}_{\beta}\setminus\{0\}\ \textrm{ for }\ i=1,\ldots, t$.
     In particular, $\NZ(\gamma;\beta)=t.$ Using inequality \ref{length},  we get
     \begin{equation}
   \begin{aligned}\label{length2}
    &\frac{\log|\gamma|}{\log|\alpha|}+e_1(\alpha)\leq m_1\leq \frac{\log|\gamma|}{\log|\alpha|}+e_2(\alpha),\\
& \frac{\log|\gamma|}{\log|\beta|}+e_1(\beta)\leq l_1\leq \frac{\log|\gamma|}{\log|\beta|}+e_2(\beta).
 \end{aligned}
 \end{equation}
Let 
    \be\label{deftheta}
      \theta= c_1\log\log|\gamma|
    \ee
    for some number $c_1$ larger than $1$ which will be chosen later. Assume that $|\gamma|>c_2$, where $c_2$ is a computable number larger than $16$.
Choose $k\in\mathbb{N}$ such that
    \be\label{defn_k}
    \theta^k\leq \frac{1}{2}\left(\frac{\log|\gamma|}{\log|\alpha|}+e_1(\alpha)\right) <\theta^{k+1}.
    \ee
Consider the intervals
    \[
    \Theta_1=(0,\theta],\ \Theta_2=(\theta,\theta^2],\ldots,\Theta_k=(\theta^{k-1},\theta^k].
    \]
We first consider the case when each of the intervals $\Theta_s$, $1\leq s\leq k$, contains at least one of the terms $m_1-m_i$, $2\leq i\leq r$ or $l_1-l_j$, $2\leq j\leq t$. Then using \eqref{deftheta} and \eqref{defn_k}, we get
    \begin{align*}
        \NZ(\gamma;\alpha)+\NZ(\gamma;\beta)&> r-1+t-1\geq k\\
        &> \frac{\log((\nicefrac{\log|\gamma|}{\log|\alpha|}+e_1(\alpha))/2)}{\log(c_1\log\log|\gamma|)}-1\\
        &\geq\frac{\log\log |\gamma|}{\log\log\log |\gamma|+C},
    \end{align*}
 for a suitable choice of $C$.
\par We next consider the case when there exists an integer $s$, $1\leq s\leq k$, such that $m_1-m_i,l_1-l_j\notin\Theta_s$ for all $i=2,\ldots,r$ and $j=2,\ldots,t$. Let $p,q$ be defined by
\begin{equation}
     \begin{aligned}\label{defn_p} 
       m_1-m_p&\leq\theta^{s-1}<\theta^s\leq m_1-m_{p+1},\\
       l_1-l_q&\leq\theta^{s-1}<\theta^s\leq l_1-l_{q+1}.
    \end{aligned}
\end{equation}  
(Here, $m_{r+1},l_{t+1}$ are defined to be zero.) 
 We have
    \be\label{split}
    (\alpha-1)\gamma=A_1\alpha^{m_p}+A_2,\ (\beta-1)\gamma=B_1\beta^{l_q}+B_2,
    \ee
where 
    \begin{align*}
        A_1&=(\alpha-1)(a_{1}\alpha^{m_1-m_p}+a_{2}\alpha^{m_{2}-m_p}+\cdots+a_{p}),\\
        A_2&=(\alpha-1)(a_{p+1}\alpha^{m_{p+1}}+\cdots+a_{r}\alpha^{m_r}),\\
       B_1&=(\beta-1)(b_{1}\beta^{l_1-l_q}+b_{2}\beta^{l_{2}-l_q}+\cdots+b_{q}),\\
        B_2&=(\beta-1)(b_{q+1}\beta^{l_{q+1}}+\cdots+b_{t}\beta^{l_t}).
    \end{align*}
    Let 
    \begin{align*}
    \epsilon_a&=|N_{K/\mathbb{Q}}(\alpha)|-1,\  \epsilon_b=|N_{K/\mathbb{Q}}(\beta)|-1,\\
        c_3&=\frac{\epsilon_a|\alpha(\alpha-1)|}{|\alpha|-1},\ c_4=|\alpha|^{-e_2(\alpha)},\\
        d_3&=\frac{\epsilon_b|\beta(\beta-1)|}{|\beta|-1},\ d_4=|\beta|^{-e_2(\beta)}.
    \end{align*}
    Now 
    \begin{equation}\label{A1bound} 
   \begin{aligned}
       &|A_1|\leq c_3|\alpha|^{m_1-m_p},\ 0<|A_2|\leq c_3|\alpha|^{m_{p+1}},\\
       &|B_1|\leq d_3|\beta|^{l_1-l_q},\ 0<|B_2|\leq d_3|\beta|^{l_{q+1}}.
    \end{aligned}      
    \end{equation} 
Further, since the $\alpha$-adic representation of the algebraic integer $A_1/(\alpha-1)$ has $m_1-m_p$ as the highest exponent, we use \eqref{length} to obtain
    \[
    |A_1|\geq c_4|\alpha-1||\alpha|^{m_1-m_p}.
    \]
    Similarly, 
    \begin{equation}
    \label{missingbound}    |B_1|\geq d_4|\beta-1||\beta|^{l_1-l_q},\ |\gamma|\geq c_4|\alpha|^{m_1},\ |\gamma|\geq d_4|\beta|^{l_1}.
    \end{equation}
Let 
    \[
    \Lambda=\frac{(\beta-1)A_1\alpha^{m_p}}{(\alpha-1)B_1\beta^{l_q}}-1.
    \]
It follows from \eqref{split} that 
    \[
    \frac{(\beta-1)(A_1\alpha^{m_p}+A_2)}{(\alpha-1)(B_1\beta^{l_q}+B_2)}=1.
    \]
We use this fact and the inequalities \eqref{A1bound}, \eqref{missingbound} and \eqref{defn_p}
to get
    \begin{align}
       \nonumber  |\Lambda|&=\left|\frac{\beta-1}{\alpha-1}\right|\left|\frac{A_1\alpha^{m_p}}{B_1\beta^{l_q}}-\frac{A_1\alpha^{m_p}+A_2}{B_1\beta^{l_q}+B_2}\right|
 =\frac{1}{|\alpha-1|}\frac{|A_1B_2\alpha^{m_p}-A_2B_1\beta^{l_q}|}{|B_1\beta^{l_q}\gamma|}\\
  \nonumber       &\leq\frac{|A_1||\alpha|^{m_p}}{|(\alpha-1)\gamma|}\cdot\frac{|B_2|}{|B_1||\beta|^{l_q}}+\frac{|A_2|}{|(\alpha-1)\gamma|}
 \\
\nonumber &\leq \frac{c_3|\alpha|^{m_1}}{|\alpha-1|c_4|\alpha|^{m_1}}\cdot\frac{d_3|\beta|^{l_{q+1}}}{d_4|\beta-1||\beta|^{l_1}}+\frac{c_3|\alpha|^{m_{p+1}}}{|\alpha-1|c_4|\alpha|^{m_1}}\\
\nonumber   &\leq \frac{c_3d_3}{|(\alpha-1)(\beta-1)|c_4d_4}\cdot\frac{1}{|\beta|^{\theta^s}}+\frac{c_3}{|\alpha-1|c_4}\cdot\frac{1}{|\alpha|^{\theta^s}}\\
\label{LFL_UB}&\leq\left(\frac{c_3d_3}{|(\alpha-1)(\beta-1)|c_4d_4}+\frac{c_3}{|\alpha-1|c_4}\right)\frac{1}{|\beta|^{\theta^s}}.
    \end{align}
    If $\Lambda\neq 0$, we apply Proposition \ref{lfl} with $T=3$, $D=2$,
    \begin{align*}
        \delta_1&=\frac{(\beta-1)A_1}{(\alpha-1)B_1}=\frac{a_{1}\alpha^{m_1-m_p}+a_{2}\alpha^{m_{2}-m_p}+\cdots+a_{p}}{b_{1}\beta^{l_1-l_q}+b_{2}\beta^{l_{2}-l_q}+\cdots+b_{q}},\ k_1=1,\\
        \delta_2&=\alpha,\ k_2=m_p,\ \delta_3=\beta,\ k_3=-l_q.
    \end{align*}
    Now, using \eqref{length2}, we have
    \begin{align}\label{LFL_B}
    \nonumber    B&=\max(|k_1|,|k_2|,|k_3|)\leq\max(m_1,l_1)\leq \frac{\log|\gamma|}{\log|\beta|}+\max(e_2(\alpha),e_2(\beta))\\
        &\leq \left(\frac{1}{\log|\beta|}+\frac{\max(e_2(\alpha),e_2(\beta)}{\log 16}\right)\log|\gamma|.
    \end{align}
    Applying inequality \eqref{height_of_product}, we obtain
    \[
        h(\delta_1)\leq h(a_{1}\alpha^{m_1-m_p}+a_{2}\alpha^{m_{2}-m_p}+\cdots+a_{p})+h(b_{1}\beta^{l_1-l_q}+b_{2}\beta^{l_{2}-l_q}+\cdots+b_{q}).
    \]
 Further, using Proposition \ref{height_of_poly} and inequality \eqref{defn_p}, we get
    \begin{align*}
        &h(a_{1}\alpha^{m_1-m_p}+a_{2}\alpha^{m_{2}-m_p}+\cdots+a_{p})\\
        &\leq (m_1-m_p)h(\alpha)+\log(|a_1|+\cdots+|a_p|)\\
     &\leq (m_1-m_p)h(\alpha)+\log(p\epsilon_a)\\
        &\leq (h(\alpha)+\log \epsilon_a+1)\max(m_1-m_p,1)
        \leq (h(\alpha)+\log \epsilon_a+1)\theta^{s-1}.
    \end{align*}
    Similarly, 
    \[
    h(b_{1}\beta^{l_1-l_q}+b_{2}\beta^{l_{2}-l_q}+\cdots+b_{q}) \leq (h(\beta)+\log \epsilon_b+1)\theta^{s-1}.
    \]
    Hence $$h(\delta_1)\leq (h(\alpha)+\log \epsilon_a+h(\beta)+\log \epsilon_b+2)\theta^{s-1}.$$ 
    Further, 
    \begin{align*}
      |\log\delta_1|&\leq|\log|\delta_1||+\pi
      \leq \log|a_{1}\alpha^{m_1-m_p}+a_{2}\alpha^{m_{2}-m_p}+\cdots+a_{p}|+\pi\\
      &\leq \log(p\,\epsilon_a|\alpha|^{m_1-m_p})+\pi\\
      &\leq(1+\log\epsilon_a+\log|\alpha|+\pi)\max(m_1-m_p,1)\leq c_6\theta^{s-1},
    \end{align*}
    where $c_6$ depends only on $\alpha$.
    Thus
    \[
    \max\left(h(\delta_1),\frac{|\log\delta_1|}{D},\frac{0.16}{D}\right)\leq c_7\theta^{s-1},
    \]
    where $c_7$ depends only on $\alpha$ and $\beta$.
     Therefore, applying Proposition \ref{lfl}, we get
    \[
    \log|\Lambda|>-c_8\theta^{s-1}\log(3eB),
    \]
    for a positive number $c_8$ depending only on $\alpha$ and $\beta$. Combining this with \eqref{LFL_UB}
 and \eqref{LFL_B}, we get that $\theta\leq c_{9}+c_{10}\log\log|\gamma|$. Since $c_{9}$ and $c_{10}$ don't depend on $c_1$, we can choose $c_1$ in equation \eqref{deftheta} so that this gives a contradiction. 
Hence $\Lambda=0$, i.e.,
    \be\label{MD1}
    \left(\frac{(\beta-1)A_1}{(\alpha-1)B_1}\right)\alpha^{m_p}\beta^{-l_q}=1.
    \ee
By Proposition \ref{multdep}, there exist integers $n_1,n_2,n_3$
such that
    \be\label{MD2}
    \left(\frac{(\beta-1)A_1}{(\alpha-1)B_1}\right)^{n_1}\alpha^{n_2}\beta^{n_3}=1
    \ee
and
    \[
    \max(|n_1|,|n_2|,|n_3|)\leq c_{11} \theta^{s-1}\leq c_{11}\theta^{k-1},
    \]
    $c_{11}$ depends only on $\alpha$ and $\beta$.
Therefore, using \eqref{defn_k}, we get
    \[
    |n_2|<\frac{1}{2}\left(\frac{\log|\gamma|}{\log|\alpha|}+e_1(\alpha)\right).
    \]
  Hence using \eqref{length2} and \eqref{defn_p}, we obtain
  \begin{align*}
      |n_2|&<\frac{1}{2}\left(\frac{\log|\gamma|}{\log|\alpha|}+e_1(\alpha)\right)\\
      &=\left(\frac{\log|\gamma|}{\log|\alpha|}+e_1(\alpha)\right)-\frac{1}{2}\left(\frac{\log|\gamma|}{\log|\alpha|}+e_1(\alpha)\right)\\
     &\leq  m_1-\theta^{s-1}\leq m_p.
  \end{align*}
Similarly, $|n_3|<l_q$. 
Using equations \eqref{MD1} and \eqref{MD2}, we get 
    \[
    \alpha^{n_2-n_1m_p}\beta^{n_3+n_1l_q}=1.
    \]
    If $n_1=0$, we get that $\alpha$ and $\beta$ are multiplicatively dependent, contradicting our assumption. If $n_1\neq 0$, then $|n_2|<m_p\leq |n_1m_p|$ and $|n_3|<l_q\leq|n_1l_q|$. Hence $n_2-n_1m_p$ and $n_3+n_1l_q$ are both non-zero. 
This again implies that $\alpha$ and $\beta$ are multiplicatively dependent, contradicting our assumption. 


\section*{Acknowledgment}
The author acknowledges the support of the DST--SERB SRG Grant (SRG/2021/000773) and the OPERA award of BITS Pilani. Thanks are due to Dr. L. Singhal for comments on an earlier draft.

\bibliographystyle{plain}
\bibliography{references}
\end{document}